\documentclass[preprint,11pt,3p]{elsarticle}
\makeatletter
\def\ps@pprintTitle{%
 \let\@oddhead\@empty
 \let\@evenhead\@empty
 \def\@oddfoot{\centerline{\thepage}}%
 \let\@evenfoot\@oddfoot}
\makeatother

\usepackage{amssymb}
\usepackage{amsmath}
\usepackage{float}
\usepackage{listings}
\usepackage{color}
\usepackage[]{algorithm2e}
\usepackage[inline]{asymptote}

\definecolor{codegreen}{rgb}{0,0.6,0}
\definecolor{codegray}{rgb}{0.5,0.5,0.5}
\definecolor{codepurple}{rgb}{0.58,0,0.82}
\definecolor{backcolour}{rgb}{0.95,0.95,0.92}
 
\lstdefinestyle{mystyle}{
    backgroundcolor=\color{backcolour},   
    commentstyle=\color{blue},
    keywordstyle=\color{red},
    numberstyle=\tiny\color{codegray},
    stringstyle=\color{codepurple},
    basicstyle=\footnotesize,
    breakatwhitespace=false,         
    breaklines=true,                 
    captionpos=b,                    
    keepspaces=true,                 
    numbers=left,                    
    numbersep=5pt,                  
    showspaces=false,                
    showstringspaces=false,
    showtabs=false,                  
    tabsize=2
}

\journal{}

\begin{document}

\begin{frontmatter}

\title{$2$-uniform words: cycle graphs, and an algorithm to verify specific word-representations of graphs}
\author[]{Ameya Daigavane, Mrityunjay Singh, Benny K. George}
\address{Indian Institute of Technology, Guwahati}
\begin{abstract}
    We prove that the number of $2$-uniform words representing the labeled $n$-vertex cycle graphs is precisely $4n$. Further, we propose a novel $O(V\log{V}+E)$-time algorithm to check whether $G(w) = G$, for a given $2$-uniform word $w$ and a graph $G = (V, E)$.  
\end{abstract}
\begin{keyword}
word-representability \sep 2-uniform word \sep cycle graph \sep fenwick tree
\end{keyword}

\end{frontmatter}
\section{Introduction}
The notion of word-representations of graphs is closely linked to the word problems of the Perkins semigroup and semi-transitive orientations of graphs, as discussed in \cite{wordsgraphs} and \cite{semitransitive}.

Consider a word $w$ on an alphabet $\Sigma$. Let the \textit{residual word} corresponding to two letters $a, b$, denoted by $w_{a, b}$, be the word obtained from the word $w$ by removing all letters of $\Sigma \setminus \{a, b\}$ from the word $w$. For example if $\Sigma=\{a,b,c,d\}$,  and $w=abaabcdbadc$ then $w_{a, c}=aaacac$. Given a word $w$, we say that letters $a$ and $b$ which appears in the word $w$ \textit{alternate} if there are no consecutive appearances of either $a$ or $b$ in the word $w_{a,b}$. For the word $w=abbcabc$, we have $w_{a,c}=acac$ implying that the letters $a$ and $c$  alternate in the word $w$. However, $w_{a,b}=abbab$, which implies that letters $a$ and $b$ do not alternate in the word $w$.
\newtheorem{dfn}{Definition}[section]
\begin{dfn}
	Let $w$ be a word on $\Sigma$. The \textit{alternating symbol graph} $G(w)=(\Sigma,E)$ is the graph in which the edge $(a,b)\in E$ iff the letters $a$ and $b$ alternate in the word $w$.
\end{dfn}

\newtheorem{example}{Example}[section]
\begin{example} Consider the word $w = bcabadc$. Its alternating symbol graph $G(w)$ is shown in Figure \ref{altgraph}.

    \begin{figure}[H]
      \label{altgraph}
      \centering
      \includegraphics[scale=0.5]{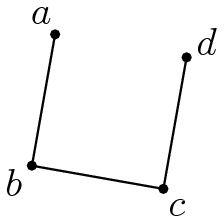}
      \caption{$G(w)$ for $w = bcabadc$.}
    \end{figure}
\end{example}

\begin{dfn}
	A graph $G=(V,E)$ is \textit{word-representable} iff there exists a word $w$ on $V$ such that $G(w) = G$. There may be many such words $w$ representing $G$. 
\end{dfn}
\begin{example}
	The first graph in Figure $2$ is word-representable, but the second graph is not, as shown in \cite{wordsgraphs}.
	
	\begin{figure}[H]
      \label{examplegraph}
      \centering
      \includegraphics[scale=0.5]{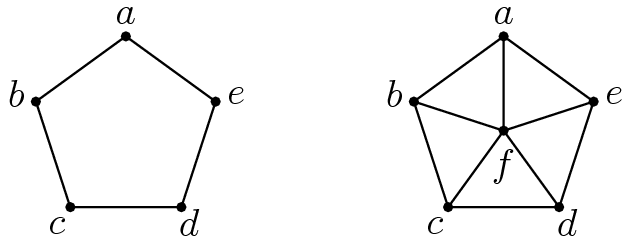}
      \caption{Word-representable and non word-representable graphs.}
    \end{figure}
    
\end{example}

Upon imposing further restrictions on the word, we get different variations of the general theme of word-representability. Here, we focus on the $2$-uniform words. A word $w$ is said to be \textit{$k$-uniform} iff every letter in $w$ occurs exactly $k$ times. If a word $w$ is $k$-uniform, then the alternating symbol graph $G(w)$ is called a $k$-word representable graph. 

Circle graphs are a special class of undirected graphs whose vertices can be associated with chords of a circle such that two vertices are adjacent iff the corresponding chords cross each other.
It has been shown in \cite{wordsgraphs} that the $2$-uniform word-representable graphs are exactly the circle graphs, excluding the complete graphs. Note that checking whether a graph with $n$ vertices is a circle graph has an $O(n^2)$-time algorithm, as described in \cite{spinrad1994recognition}. 

It is simple to see that the cycle graphs $C_n$ are a subset of the circle graphs. An algorithm has been given in \cite{wordsgraphs} to construct a word $w_n$ representing a labeled cycle graph on $n$ vertices. We build upon this in section $2$ and show that exactly $4n$ $2$-uniform words represent a labeled cycle graph, all obtainable from cyclic shifts and reflection of the word $w_n$. 

In section $3$, we describe and prove the $O(V \log{V} +E)$-time algorithm employing Fenwick Trees to check if a $2$-uniform word $w$ satisfies $G(w) = G$ for a given graph $G=(V,E)$. 

\newtheorem{thm}{Theorem}
\newproof{pot}{Proof}
\section{$2$-uniform Representations of the Cycle Graphs}
Before we investigate the cycle graphs, let us show how the $2$-uniform words are related to the circle graphs described above. We can represent a $2$-uniform word $w$ of length $l$ on a circle, labelled by positions from $1$ to $l$, clockwise. At each position, we mark the letter in $w$ that occurs at this position. Now, join the two points where a specific letter exists by a chord. Note that two letters alternate iff their corresponding chords intersect.
We call this the \textit{circle representation} of $w$.

Consider the cycle graph $C_n$ labelled $1, 2...n$ in the clockwise direction, where $n > 3$. Note that vertices $k$ and $k + 1$ are connected for $1 \leq k \leq n - 1$, and the edge from $n$ to $1$ completes the cycle. 
Define the word $w_n$, on the alphabet $\mathcal{A} = \{1, 2, .. n\}$, by
\[ w_n = 1 n 2 1 3 2 4 3 5 4 ... (n - 1) (n - 2) n (n - 1). \]
It can be easily verified that $G(w_n) = C_n$. (Between every two occurrences of letter $r$, we have exactly $(r - 1)$ and $(r + 1)$ where addition and subtraction are done cyclically in $[n]$.)

Hence, every cycle graph has a $2$-uniform word representation. We claim that $w_n$, in fact, 'generates' all the other $2$-uniform word representations of $C_n$.

\begin{figure}[H]
      \centering
      \includegraphics[scale=0.5]{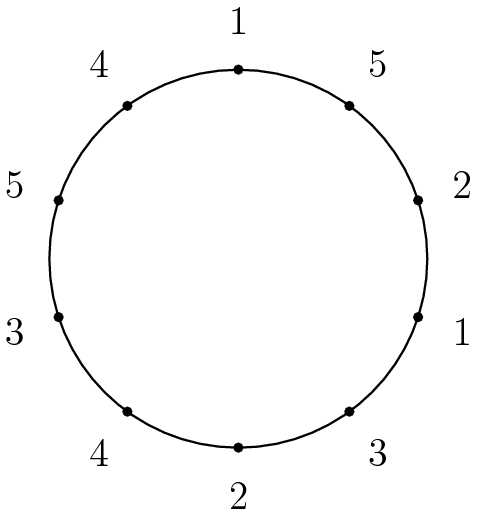}
      \caption{The circle representation of $w_5$. We can start at any of the $10$ positions, and read in either the clockwise or the anticlockwise direction, to get a total of $4 \times 5 = 20$ distinct words.}
\end{figure}

\newtheorem{lem}{Lemma}
\begin{lem}
Given a $k$-uniform word $w$ that represents a graph $G$, the word $w^{'}$ given by any cyclic shift (a 'rotation'), or a reflection of $w$, also represents $G$.
\end{lem}

\begin{pot}
Under rotations, it is not too difficult to check that an alternation of $c_1$ and $c_2$ in $w$ is preserved in $w^{'}$, and no new alternations are created in $w^{'}$.
Similarly, under a reflection, every alternation of $c_1$ and $c_2$ is retained (the residual word $w_{c_1, c_2}$ has the letters $c_1, c_2$ simply interchanged with no new alternations created.)
As the alternations are the same in $w$ and $w^{'}$, the graphs they represent are the same too.
\end{pot}

\begin{lem}
Each composition of one of the $2n$ rotations and a reflection of $w_n$ gives rise to a distinct word. 
\end{lem}

\begin{pot}
We provide an outline of the proof. Clearly, each rotation of $w_n$ gives rise to a distinct word. (Look at the starting letter, say $l$. If those are not different, then look at the immediate right letter. In one word, this will be $l - 1$, and the other, this will be $l + 2$ - both addition and subtraction done cyclically. These are distinct, because $n > 3$.)

Note that $w_n$ has the property that two occurrences of the letter $r$ are separated by exactly 2 letters, $(r + 1)$ and $(r - 1)$ in exactly this order (if counting cyclically, in the forward direction.) Reflections and rotations do not alter this distance, however, reflections flip the order of $(r + 1)$ and $(r - 1)$, while rotations maintain this. As $r$ is arbitrary, we have that a rotation followed by a reflection of $w$ cannot be a rotation of $w$.
\end{pot}

\begin{lem}
Let $w$ be any $2$-uniform word that represents $C_n$ labelled $1, 2, 3 ... n$ for $n > 3$. Further, let $w[i]$ be the letter at position $i$ for all $ 0 \leq i \leq 2n - 1$. 
\noindent Then, the circle representation of $w$ satisfies the following property:
\noindent For every $r \in \{1, 2, 3 ... n\}$, the two sets of positions,
\[ U_r = \{i : (w[i] - r) > 1\}  \text{ and }  L_r = \{ i : (r - w[i]) > 1\},\] 
if both are non-empty, lie entirely in one of the two segments defined by the chord corresponding to $r$. If exactly one is non-empty, then that set lies entirely in one of two segments. 
\begin{figure}[H]
  \centering    
  \includegraphics[scale=0.4]{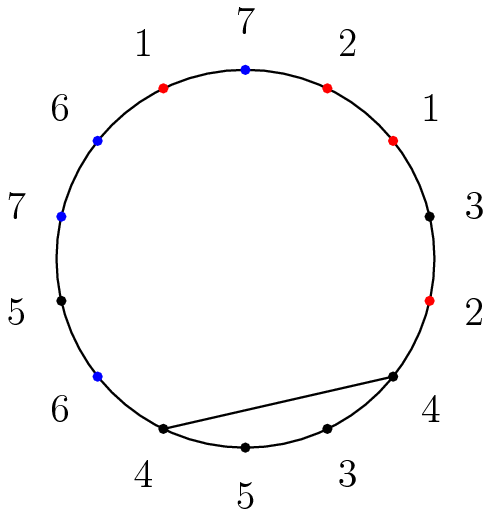}
  \caption{Circle representation of 12132546576734, representing $C_7$. The chord corresponding to $r = 4$ has been drawn. The two sets of points have been coloured in red and blue, respectively. Black points belong to neither of the two sets.}
\end{figure}
\end{lem}
\begin{pot}
Consider an arbitrary $r \in \{1, 2, 3 ... n\}$. Note that exactly one of two sets described is empty when $r = 1, 2, n - 1$ or $n$. When $3 \leq r \leq n - 2$, then both sets are non-empty. We consider this case first.
Call the segments defined by the chord corresponding to $r$ as the 'upper segment' and 'lower segment', as shown in the figure.
There must be an $r + 1$ and an $r - 1$ on opposite sides of the chord corresponding to $r$.
Without loss of generality, we can assume the $r + 1$ before the $r - 1$, when traversing along the circle, clockwise from the left $r$ to the right one. (This is because $l$ under reflection represents the same graph.)
Note that this fixes the positions of $r + 1$ and $r - 1$ on the other side of the chord. (As $n > 3$, $r + 1$ and $r - 1$ cannot alternate.)

\begin{figure}[H]
  \includegraphics[scale=0.4]{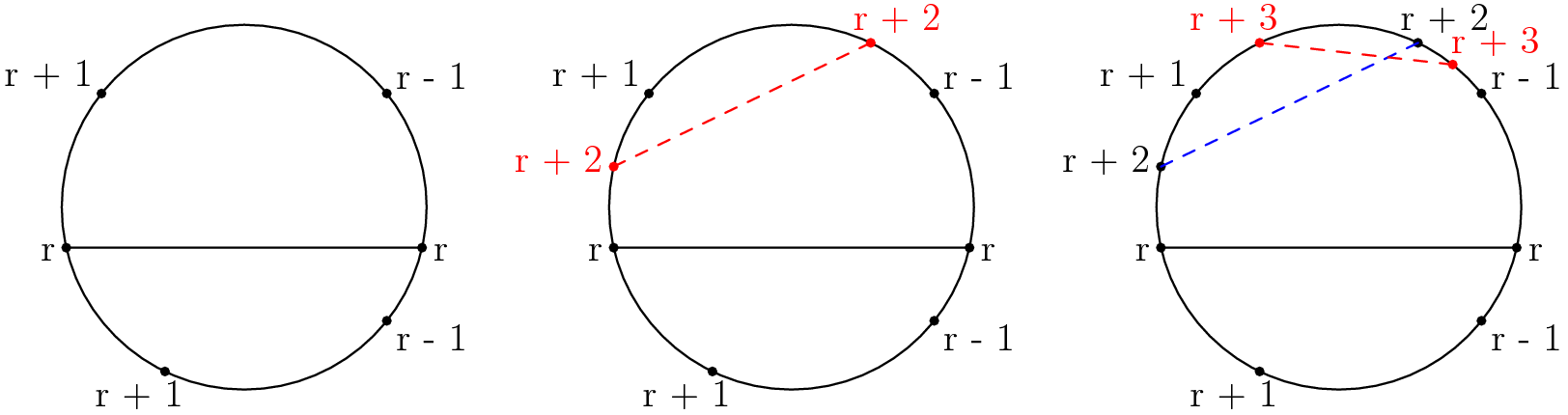}
  \caption{The steps of the proof shown visually.}
\end{figure}
Now, $r + 1$ must alternate with $r + 2$. This means there must be an $r + 2$ in between the two $r + 1$s on either side.
But $r + 2$ cannot alternate with $r$ (again, as $n > 3$). This forces the chord corresponding to $r + 2$ to lie completely in the upper segment, or completely in the lower segment (depending on where we keep the left $r + 2$). The 'upper segment' choice is shown in Figure 5 above.

Similarly, if $r \leq n - 3$, an $r + 3$ must lie in between the two $r + 2$s, and again cannot alternate with $r$.
Thus, the chords corresponding to $r + 2$ and $r + 3$ lie in the upper segment.
Extending the above to $r + 4, r + 5..$, it is clear that for every $s$, $r + 2 \leq s \leq n$, $s$ will lie in the upper segment. But since $s$ cannot alternate with $r$, the chord corresponding to $s$ lies completely in the upper segment.
Therefore, $U_r$ lies completely in the upper segment. (If we had chosen the first $r + 2$ to lie in the lower segment, then $U_r$ would lie completely in the lower segment.)

Now, we can place $r - 2$ in the upper segment or lower segment, in between the two $r - 1$s ($r - 2$ again cannot alternate with $r$). Again,  for every $s$, $ 1 \leq s \leq r - 2$, $s$ will lie in the same segment as the $r - 2$s. We conclude that $L_r$ too lies completely either in the upper or the lower segment.

Suppose $L_r$ lies in a different segment than $U_r$. As $L_r$ definitely contains the two $1$s, and $U_r$ definitely contains the two $n$s, the two $1$s and two $n$s lie in distinct segments. This means they cannot alternate, and hence the undirected edge $(1, n)$ would not be present in $G(w)$, making it impossible for $w$ to represent $C_n$.

If one of $U_r$ and $L_r$ is empty, then our reasoning above shows that this set must lie entirely in one of the segments.

Thus, if both are non-empty, $U_r$ and $L_r$ must lie entirely in the same segment, and if one is non-empty, then that set lies entirely in one segment - proving our claim.
\end{pot}

\begin{thm}
The cycle-graph $C_n$ (where $n > 3$) has exactly $4n$ $2$-uniform word representations, each given by a rotation, or a rotation followed by a reflection, of $w_n$.
\end{thm}

\begin{pot}
As $w_n$ is $2$-uniform, any rotation, or any rotation followed by a reflection of $w_n$, generates $C_n$ as well, by Lemma 1. There are $2n$ possible rotations (including the 'zero' rotation that fixes $w_n$) and a possible reflection, each giving $2n \times 2 = 4n$ words, each distinct by Lemma 2.

We show that these $4n$ words are the only possibilities. 
Consider any $2$-uniform word $w$ representing $C_n$. Take an arbitrary $r \in \{1, 2, 3 ... n\}$, and look at the two positions of it in the circular representation of $w$. As defined in Lemma 3, $L_r$ and $U_r$ both lie in one of the segments, on one side of the chord corresponding to $r$. Call this segment, Segment 1. The only points that are not included in $L_r \cup U_r$, are the $2$ letters to which $r$ is connected - and these must be included in Segment 2. (For example, if $n = 10$, for $r = 5$, these are $4$ and $6$, for $r = 1$, these are $10$ and $2$.)
Thus, in the circle representation of $w$, between two positions of $r$, there are exactly two elements in one direction - exactly the elements to which $r$ is connected to in $C_n$.

Fixing the positions of $1$ on the circle representation of $w$, we have either $2$ then $n$ (giving $w = ..12n1..$), or $n$ then $2$ (giving $w = ..1n21..$), between them. 
If we take $..12n1..$, then the position right after the right $1$ cannot be a $2$, as between two $2$s there must be a $3$. Therefore, $w = ..2312n1..$. Now, as there is no $3$ three places to the right of the current $3$, $3$ must be present three places to the left of the current $3$. As there must be a $4$ in between two $3$s, we must have $w = ..342312n1...$.
We can continue this way, until all letters are accounted for.
Similarly, if we take $..1n21..$, then we can continue filling in positions, with the same reasoning as above.

Note that in either case, the word $w$ obtained is in fact, a rotation, or a rotation after a reflection of $w_n = 1 n 2 1 3 2 4 3 5 4 ... (n - 1) (n - 2) n (n - 1).$
(The starting positions of $1$ on the circle representations can be shifted, and these all represent the rotations of $w_n$.)
Thus, every word representing $C_n$ must be one of the $4n$ distinct $2$-uniform word representations, obtained from $w_n$, as claimed.
\end{pot}

\section{Is $G(w) = G$?}
Given a labelled graph $G=(V,E)$, and a $2$-uniform word $w$ on $V$, we can ask if $w$ is a word-representation of $G$; essentially, is $G(w) = G$? Below is our proposed $O(V \log{V} + E)$-time algorithm employing Fenwick Trees (described in \cite{fenwick}), with $O(V)$ auxiliary space.

\begin{algorithm}[H]
\DontPrintSemicolon
 \KwData{$2$-uniform word $w$ on $V$, and  graph $G=(V,E)$.}
 \KwResult{Returns \textbf{true} if $G(w) = G$, and \textbf{false} otherwise.}
 \;
 Initialize FenwickTree with 0 in all positions with total length $w.length()$. \;
 Initialize array of positions $pos[]$ to (NULL, NULL) for all letters in $w$.\;
 edgecount = 0\;
 \;
 \For{k $=0$ to $w$.length() $- 1$}{
  \eIf{$pos[w[k]].first = NULL$}{
       $pos[w[k]].first = k$
       \tcp*{$w[k]$ appears for the first time}
   }{
    $pos[w[k]].second = k$
    \tcp*{$w[k]$ appears for the second time}
    $i = pos[w[k]].first$\;
    $j = pos[w[k]].second$\;
    \;
    \tcp{add the number of unmarked nodes in $w[i...j]$}
    edgecount += $j - i - \text{FenwickTree.rangesum}(i + 1, j - 1) - 1$ 
    \;
    \;
    \tcp{mark the positions $i$ and $j$}
    FenwickTree.update($i, 1$) \;
    FenwickTree.update($j, 1$) \;
  }
 }
 
 \eIf{edgecount $\neq |E|$}{
 return false
 }{
 \For{edge $(u, v)$ in $E$}{
    \If{u and v do not alternate}
    {   
    return false \tcp*{only a $O(1)$ comparison using $pos[u]$ and $pos[v]$}
    }
 }
 return true

 }
 \caption{GraphCheck}
\end{algorithm}

\subsection*{Time and Space Complexity:}
\noindent We mark the positions of letters already considered.
Using Fenwick Trees, obtaining the number of marked nodes in a range and updating the marked nodes both take $O(\log{V})$ time. As we have to do this $V$ times, once for each letter, but scanning the word only once, the overall time complexity of the first for-loop is $O(V\log{V})$.
The second for-loop takes $O(E)$ time, as it only takes constant time to check if a certain edge exists, with the array of positions, $pos[]$.
Thus, the overall time complexity of this algorithm is $O(V\log{V} + E)$, as claimed.

Note that the Fenwick Tree and the array of positions can be implemented in $O(V)$ space. Thus, our algorithm requires $O(V)$ auxiliary space. 

Note that naively checking if every edge in $G$ in present in $G(w)$ would be an $O(VE)$-time algorithm - for each edge $(u, v)$ in $G$, scan the word $w$ to check if $u$ and $v$ alternate. Hence, our algorithm is an aymptotic improvement.

\begin{pot}
A letter $r$ alternates with another letter, say $l$, iff there is exactly one occurrence of $l$ in between the two occurrences of $r$. 
Thus, we scan the word $w$, finding the first letter occurring twice. This letter must have edges with all letters occurring in between its two occurrences, so we add the number of 'unmarked' letters between them to the edgecount of $G(w)$. As this letter cannot have any more edges incident to it, we 'mark' the positions of this letter, updating the Fenwick tree, and continue until the entire word is scanned.
We first check that the number of edges in $G$ and $G(w)$ are the same.
Then, we check if every edge in $G$ is in $G(w)$. Together, this proves that $G(w) = G$.
\end{pot}

An implementation of the algorithm above, in Python 3.5, is provided in the reference \cite{pycode} below.

\bibliographystyle{elsarticle-num}

\bibliography{sample}

\end{document}